\documentclass{article}
\usepackage{amssymb,amsmath}
\usepackage{graphicx}

\def\qed{\hfill {\hbox{${\vcenter{\vbox{               
   \hrule height 0.4pt\hbox{\vrule width 0.4pt height 6pt
   \kern5pt\vrule width 0.4pt}\hrule height 0.4pt}}}$}}}

\newtheorem{theorem}{Theorem}
\newtheorem{definition}{Definition}

\newtheorem{example}{Example}

\newtheorem{remark}{Remark}

{\qed\par\medskip}

\date{}

\title{\Large \textbf{Polynomial Birack Modules}}

\author{Evan Cody\footnote{Email: evanwcody@gmail.com}\and
Sam Nelson\footnote{Email: knots@esotericka.org}}

\begin{document}
\maketitle

\begin{abstract}
Birack modules are modules over an algebra $\mathbb{Z}[X]$ 
associated to a finite birack $X$. In previous work, birack module 
structures on $\mathbb{Z}_n$ were used to enhance the birack counting
invariant. In this paper, we use birack modules over Laurent 
polynomial rings $\mathbb{Z}_n[q^{\pm 1}]$ to enhance the birack counting
invariant, defining a customized Alexander polynomial-style signature for
each $X$-labeled diagram; the multiset of these polynomials is an enhancement
of the birack counting invariant. We provide examples to demonstrate that 
the new invariant is stronger than the unenhanced birack counting invariant
and is not determined by the generalized Alexander polynomial.
\end{abstract}

\medskip

\quad
\parbox{5in}{
\smallskip

}

\section{\large\textbf{Introduction}}\label{I}

\textit{Biracks} are algebraic structures with axioms motivated by framed 
oriented Reidemeister moves. They were introduced in \cite{FRS} in order to 
define invariants of framed oriented knots and links. In \cite{N2} a property 
of finite biracks called \textit{birack rank} was used to define a computable 
integer-valued invariant of unframed oriented knots and links called the 
\textit{integral birack counting invariant}, $\Phi_X^{\mathbb{Z}}$.
In \cite{AG} an algebra $\mathbb{Z}[X]$ called the \textit{rack algebra} was 
associated to a finite rack $X$ (a particular type of birack), and modules 
over $\mathbb{Z}[X]$ were studied. In \cite{HHNYZ} rack module structures
over $\mathbb{Z}_n$ were used to enhance the rack counting invariant, 
defining a new invariant $\Phi_X^M$ which specializes to the integral rack 
counting invariant $\Phi_X^{\mathbb{Z}}$ but is generally stronger. In \cite{BN} 
rack modules were generalized to the case of biracks, and birack modules 
over $\mathbb{Z}_n$ were employed to define enhancements of the birack 
counting invariant.

In this paper we consider birack modules over Laurent polynomial rings 
$\mathbb{Z}_n[q^{\pm 1}]$; such a birack module lets us define a customized
Alexander polynomial for each birack homomorphism $f:BR(L)\to X$, the multiset 
of which forms a new enhancement of the integral counting invariant.
Moreover, the generalized  Alexander (Sawollek) and classical Alexander 
polynomials emerge as special cases of the enhanced invariant.

The paper is organized as follows. In Section \ref{B} we review the basics
of biracks and birack modules. In Section \ref{DC} we define the enhanced
link invariant. In Section \ref{CA} we collect some examples 
and in Section \ref{Q} we end with some questions for future work.

\section{\large\textbf{Biracks and birack modules}}\label{B}

We begin with a definition. First introduced in \cite{FRS}, a \textit{birack}
is an algebraic structure consisting of a set $X$ and a map 
$B:X\times X\to X\times X$ with axioms derived from the \textit{oriented
framed Reidemeister moves}, obtained by considering all ways of orienting
the strands in the moves
\[
\includegraphics{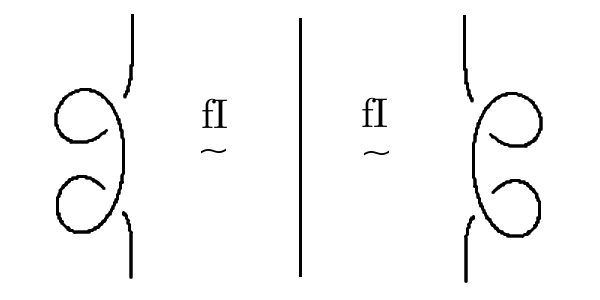}\quad
\includegraphics{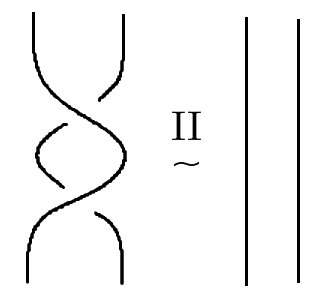}\quad
\includegraphics{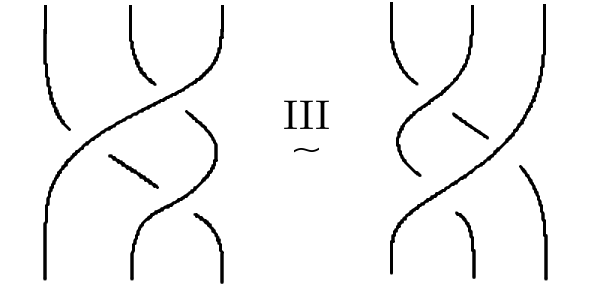}
\]
with the correspondence (also known as the \textit{semiarc labeling rule})
\[\raisebox{-0.5in}{\includegraphics{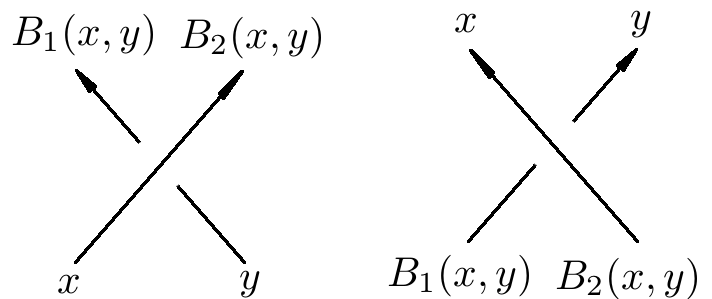}} 
\quad B(x,y)=(B_1(x,y),B_2(x,y))\]

In \cite{KR} and later in \cite{FJK} the unframed oriented case, known as
the \textit{strong biquandle} case, was
considered; the version below comes from \cite{N2}.

\begin{definition}
\textup{Let $X$ be a set and $\Delta:X\to X\times X$ the diagonal map defined
by $\Delta(x)=(x,x)$. Then an invertible map $B:X\times X\to X\times X$ is a 
\textit{birack map} if the following conditions are satisfied:
\begin{itemize}
\item[(i)] There exists a unique invertible map $S:X\times X\to X\times X$
called the \textit{sideways map} such that for all $x,y\in X$, we have
\[S(B_1(x,y),x)=(B_2(x,y),y).\]
\item[(ii)] The components $(S^{\pm 1}\Delta)_{1,2}:X\to X$ of the composition 
of the sideways map and its inverse with the diagonal map $\Delta$ are 
bijections, and
\item[(iii)] $B$ satisfies the \textit{set-theoretic Yang-Baxter equation}
\[(B\times I)(I\times B)(B\times I)=(I\times B)(B\times I)(I\times B).
\]
\end{itemize}
}\end{definition}

These axioms are the conditions required to ensure that each labeling of an 
oriented framed knot or link diagram according to the semiarc labeling rule 
above before a framed oriented Reidemeister move corresponds to a unique 
such labeling after the move. By construction, the number of labelings of 
an oriented framed link diagram is an invariant of framed isotopy, called the
\textit{basic birack counting invariant}, denoted $\Phi_X^{B}$.
 
Let $X$ be a birack. The bijection $\pi:X\to X$ defined by 
$\pi=(S\Delta)_1(S\Delta)_2^{-1}$ represents going through a positive kink:
\[\includegraphics{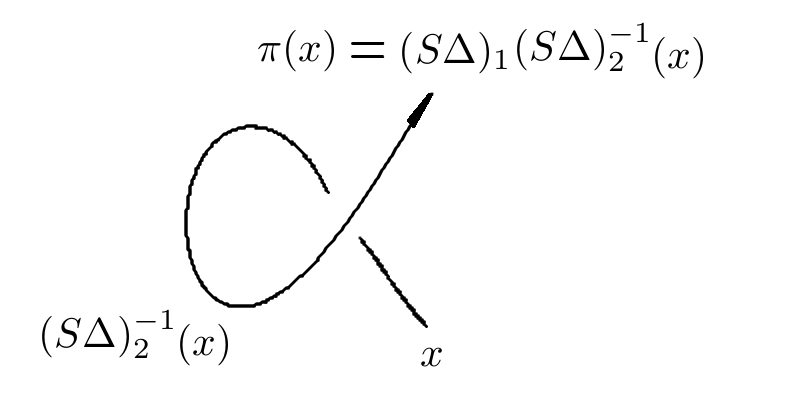}\]
This bijection $\pi$, called the \textit{kink map}, 
is an element of the symmetric group on $X$; if $X$ is finite, then $\pi$
has a finite exponent $N\in \mathbb{Z}$ such that $\pi^N(x)=x$ for all 
$x\in X$. We call $N$ the \textit{birack rank} or \textit{birack 
characteristic} of $X$. 
The map $(S\Delta)_2^{-1}$ is sometimes called $\alpha$.
A birack with rank $N=1$ is a \textit{strong biquandle}.

Recall that the framing number of a component of a blackboard framed link
is given by the \textit{writhe} of the component, i.e., the number of positive 
self-crossings of the component minus the number of negative self-crossings. 
\[\includegraphics{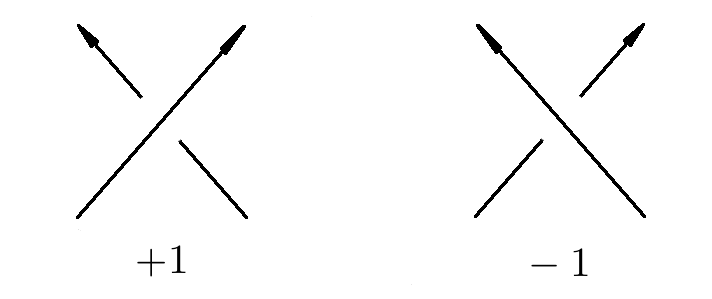}\]
On a link of $c$ components, we can specify framings with framing vectors
$\vec{w}\in\mathbb{Z}^c$ where the $k$th component of $\vec{w}$ is the number 
of times the $k$th component crosses over itself at a positive crossing
minus the number of times the $k$th component crosses over itself at a 
negative crossing.

If $X$ is a finite birack with rank $N$, then labelings of a link by $X$
before and after an \textit{$N$-phone cord move}
\[\includegraphics{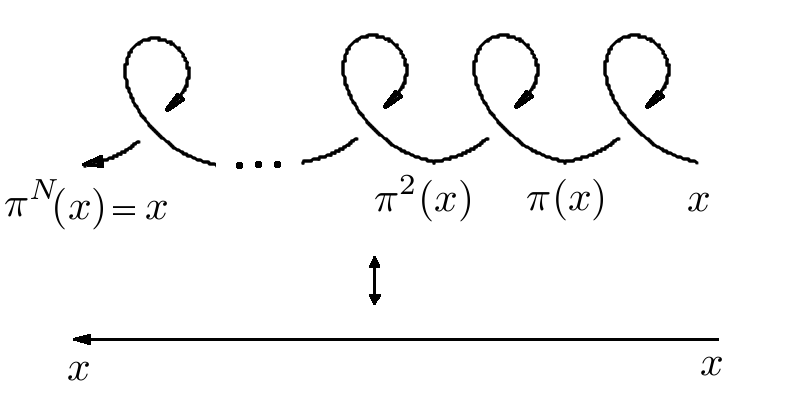}\]
are in bijective correspondence. 
Thus, for any oriented link $L$ of $c$ components, the $c$-dimensional 
lattice of basic counting invariant values $\Phi_X^B(L,\vec{w})$ is tiled 
with a tile consisting of framing vectors in $\{0,1,2,\dots,n-1\}^c$ (which
we may identify with
$(\mathbb{Z}_N)^c$). The sum of 
basic counting invariants over one such tile is an invariant of the unframed 
oriented link known as the \textit{integral birack counting invariant},
\[\Phi_X^{\mathbb{Z}}=\sum_{\vec{w}\in(\mathbb{Z}_N)^c} \Phi_X^B(L,\vec{w}).\] 
In particular, we have

\begin{theorem} If $X$ is a finite birack and $L$, $L'$ are ambient isotopic
oriented links, then $\Phi_X^{\mathbb{Z}}(L)=\Phi_X^{\mathbb{Z}}(L')$.
\end{theorem}
See \cite{N2} for more.

We can specify a birack map on a finite set $X=\{x_1,x_2,\dots,x_n\}$ with a
\textit{birack matrix} specifying the operation tables of the two component
maps $B_1,B_2:X\times X\to X$ of $B$ considered as binary operations. More 
precisely, the birack matrix of $X$ is an $n\times 2n$ block matrix $[U|L]$ 
with $U_{j,k}=l$ and $L_{j,k}=m$ where $x_l=B_1(x_k,x_j)$ and $x_m=B_2(x_j,x_k)$. 
Note the reversed order of the inputs in $B_1$; the notation is chosen so 
that the output label and the input row label are part of the same strand. 
It is sometimes convenient to abbreviate $B_1(x,y)=y^x$ and $B_2(x,y)=x_y$,
i.e. $B(x,y)=(y^x,x_y)$.

\begin{example}\label{ex-1}
\textup{Let $X=\{x_1,x_2\}$ be the birack with birack matrix 
\[M=\left[\begin{array}{rr|rr}
1 & 1 & 2 & 2 \\
2 & 2 & 1 & 1 \\
\end{array}\right].\] We can interpret $X$ in this example as a labeling rule
saying that when a strand crosses under another strand, it keeps the
same label, but a strand crossing over another strand switches labels
from 1 to 2 or from 2 to 1. The kink map $\pi$ is then the transposition $(12)$,
so this birack has rank $N=2$.}

\textup{To compute the integral birack counting invariant $\Phi_X^{\mathbb{Z}}$ 
for a link $L$ and birack $X$ of rank $N$, we need to consider diagrams of 
$L$ with all framing vectors in $(\mathbb{Z}_N)^c$. For a knot $K$ and our 
birack $X$, this means finding all labelings of one diagram of $K$ with an 
even writhe and one diagram of $K$ with an odd writhe. For example, the 
even-writhe figure eight knot has two 
$X$-labelings, while the odd-writhe figure eight knot has no valid 
$X$-labelings.
\[\includegraphics{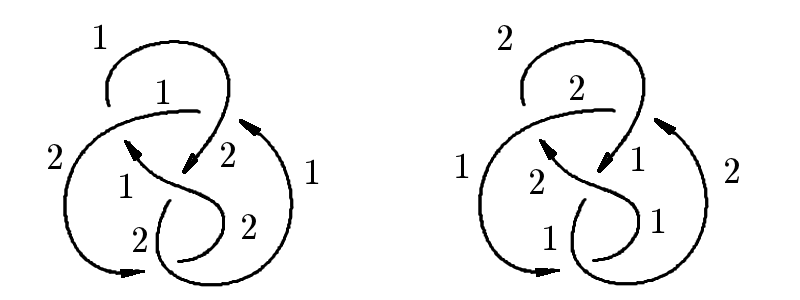} \]
Similarly, the even-writhe unknot has two $X$-labelings and the odd-writhe 
unknot has no $X$-labelings.}
\[\includegraphics{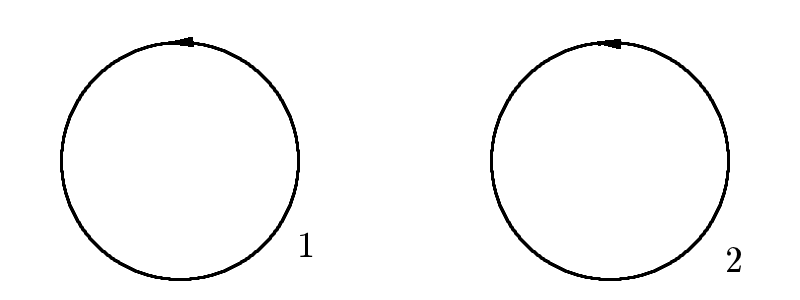}\]
\textup{For a link of two components, we would need to count
labelings for a diagram of $L$ with both writhes even, one with both writhe odd,
and the two even-odd writhe combinations, etc.}
\end{example}

The integral birack counting invariant $\Phi^{\mathbb{Z}}_X$ in example \ref{ex-1} 
does not distinguish the figure eight knot from the unknot, but these 
$X$-labeled knot diagrams are still apparently quite different. Thus, we would 
like to develop a way to tell the birack-labeled links apart. An invariant
of birack-labeled links is called a \textit{signature} of the
labeled link; collecting the signatures of the labelings of a link
over a complete $(\mathbb{Z}_N)^c$ tile defines an \textit{enhancement}
of the integral birack counting invariant. 
In \cite{BN} an enhancement of the birack counting invariant was defined using
\textit{birack modules}:

\begin{definition}\label{def2}\textup{
Let $X$ be a finite birack and let $\Lambda$ be the
polynomial ring $\Lambda=\mathbb{Z}[t_{x,y}^{\pm 1}, s_{x,y},r_{x,y}^{\pm 1}]$
with commuting invertible variables $t_{x,y},r_{x,y}$ and generic variables 
$s_{x,y}$ indexed by ordered pairs of birack elements $x,y\in X$. The 
\textit{birack algebra} associated to $X$ is the quotient algebra 
$\mathbb{Z}[X]=\Lambda/I$ where $I$ is the ideal generated by elements of 
the forms}
\[
\begin{array}{rlrlrl}
\bullet & r_{x_y,z}r_{x,y}-r_{x_{z^y},y_z}r_{x,z^y} &
\bullet & t_{x_{z^y},y_z}r_{y,z}-r_{y^x,z^{x_y}}t_{x,y} &
\bullet & s_{x_{z^y},y_z}r_{x,z^y}-r_{y^x,z^{x_y}}s_{x,y} \\
\bullet & t_{x,z^y}t_{y,z}-t_{y^x,z^{x_y}}t_{x_y,z} &
\bullet & t_{x,z^y}s_{y,z}-s_{y^x,z^{x_y}}t_{x,y} &
\bullet & s_{x,z^y}-t_{y^x,z^{x_y}}s_{x_y,z}r_{x,y}-s_{y^x,z^{x_y}}s_{x,y} \\
\multicolumn{6}{l}{\bullet \quad  1-\displaystyle{\prod_{k=0}^{N-1} (t_{\pi^k(x),\alpha(\pi^k(x))}r_{\pi^k(x),\alpha(\pi^k(x))}+s_{\pi^k(x),\alpha(\pi^k(x))})}} 
\end{array}\]
\noindent
\textup{for $x,y,z\in X$ where we have $B(x,y)=(y^x,x_y)$.
}\end{definition}

If $\mathbf{R}$ is a commutative ring, we can give $\mathbf{R}$ the structure 
of a $\mathbb{Z}[X]$-module by choosing elements $t_{x,y},s_{x,y},r_{x,y}$
of $\mathbf{R}$ such that the ideal $I$ in $\mathbf{R}$ is zero. We can specify 
such a structure with an $n\times 3n$ block matrix $M_{\mathbf{R}}=[T|S|R]$ whose
entries $T_{i,j}=t_{x_i,x_j},S_{i,j}=s_{x_i,x_j}$ and $R_{i,j}=r_{x_i,x_j}$ make each of
the generators of $I$ zero.

\begin{example}\label{ex0}\textup{
The birack $X$ with matrix $M_X=\left[\begin{array}{rr|rr}
1 & 1 & 2 & 2 \\
2 & 2 & 1 & 1
\end{array}\right]$ has $\mathbb{Z}[X]$-modules over the ring
$\mathbf{R}=\mathbb{Z}_5$ including the module $M$ given by the matrix
\[M_{\mathbf{R}}=\left[\begin{array}{rr|rr|rr}
1 & 1 & 2 & 1 & 2 & 2 \\
1 & 1 & 4 & 2 & 3 & 3 \\
\end{array}\right].\]}
\end{example}

Given an $X$-labeling $f$ of an oriented framed link diagram $L$ and a 
$\mathbb{Z}[X]$-module $M$, we can give the semiarcs in $L$ a secondary 
labeling by elements of $M$, usually visualized as beads:
\[\includegraphics{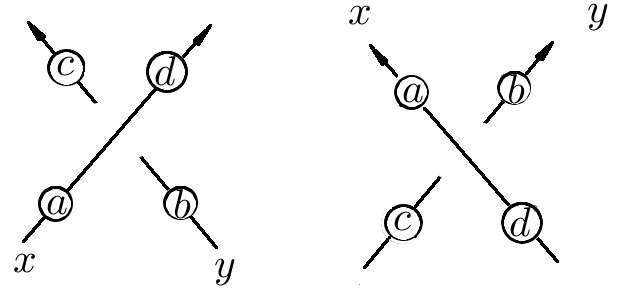}\quad \quad \quad
\raisebox{0.5in}{$\begin{array}{rcl}
c& = & t_{x,y} b+s_{x,y} a \\
d & = & r_{x,y} a. \\
\end{array}$}\]
The ideal $I$ in definition \ref{def2} is chosen so that bead labelings before and after $X$-labeled
framed Reidemeister and $N$-phone cord moves are in one-to-one correspondence.
For each $X$-labeling, the set of bead labelings by $M$ forms a 
$\mathbb{Z}[X]$-module with presentation matrix determined by the crossing 
relations, called the \textit{fundamental $\mathbb{Z}[X]$-module} of $f$,
denoted $\mathbb{Z}[f]$. 
Replacing the variables $t_{x,y}, s_{x,y},r_{x,y}$ with their values in $M$
yields the \textit{fundamental $M$-module} of the birack labeled diagram $L$,
denoted $M[f].$ By construction, $X$-labeled framed Reidemeister and $N$-phone
cord moves induce isomorphisms of $M[f]$.

\begin{example}\label{ex2}\textup{
The $X$-labeled figure eight knot below has fundamental $\mathbb{Z}[X]$-module
and fundamental $M$-module given by the listed presentation matrices where
the semiarcs are numbered starting with the pictured basepoint following the
orientation. Replacing the variables $t_{x,y}$, $s_{x,y}$ and $r_{x,y}$ with their
values in the $X$-module in example \ref{ex0} yields the following matrix with
entries in $M=\mathbb{Z}_5$:
\[\raisebox{-0.5in}{\includegraphics{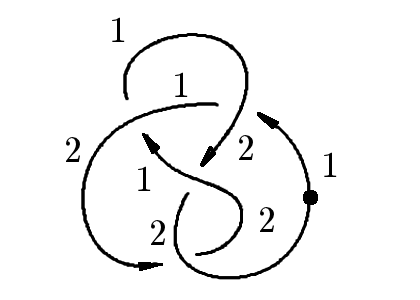}}
\begin{array}{cc}
\mathbb{Z}[f] & M[f] \\ \hline 
&  \\
\left[\begin{array}{rrrrrrrr}
-1 & t_{2,1} & 0 & 0 & 0 & 0 & s_{2,1} & 0 \\
0 & 0 & 0 & 0 & 0 & -1 & r_{2,1} & 0 \\
0 & 0 & s_{2,1} & 0 & -1 & t_{2,1} & 0 & 0 \\
0 & -1 & r_{2,1} & 0 & 0 & 0 & 0 & 0 \\
0 & 0 & 0 & s_{2,2} & 0 & 0 & t_{2,2} & -1 \\
0 & 0 & 0 & r_{2,2} & -1 & 0 & 0 & 0 \\
0 & 0 & t_{2,2} & -1 & 0 & 0 & 0 & s_{2,2} \\
-1 & 0 & 0 & 0 & 0 & 0 & 0 & r_{2,2}
\end{array}\right] &
\left[\begin{array}{rrrrrrrr}
4 & 1 & 0 & 0 & 0 & 0 & 4 & 0 \\
0 & 0 & 0 & 0 & 0 & 4 & 3 & 0 \\
0 & 0 & 4 & 0 & 4 & 1 & 0 & 0 \\
0 & 4 & 3 & 0 & 0 & 0 & 0 & 0 \\
0 & 0 & 0 & 2 & 0 & 0 & 1 & 4 \\
0 & 0 & 0 & 3 & 4 & 0 & 0 & 0 \\
0 & 0 & 1 & 4 & 0 & 0 & 0 & 2 \\
4 & 0 & 0 & 0 & 0 & 0 & 0 & 3
\end{array}\right]
\end{array}
\]
Row-reduction over $\mathbb{Z}_5$ yields a 3-dimensional solution space,
so there are $5^3=25$ total $M$-labelings of the pictured $X$-labeled
diagram. 
}\end{example}

In \cite{BN}, the \textit{birack module enhanced invariant} associated to
a pair $(X,M)$ of a birack $X$ and birack module $M$ was defined as
\[\Phi_{X}^{M}(L)=\sum_{f\in \mathcal{L}(L,X)} u^{|M[f]|}\]
where $\mathcal{L}(L,X)$ is a set of $X$-labelings over a complete tile
of framings of $L$ modulo the rank rank $N$ of $X$.

\begin{example}\label{ex3}\textup{
The birack module enhanced invariant associated to the birack and module 
in example \ref{ex2} distinguishes the figure eight knot $4_1$ from the 
unknot $0_1$ with $\Phi_X^M(4_1)=2u^{25}\ne 2u^5=\Phi_X^M(0_1)$. 
}\end{example}

\section{\large\textbf{Polynomial birack modules and an enhanced
link invariant}}\label{DC}

Let $\mathbf{R}$ be a Laurent polynomial ring and let $M$ be a module over
$\mathbf{R}$ with presentation matrix $A$, i.e., 
$A\in M_{n,m}(\mathbf{R})$ such that $M=\mathrm{CoKer}(A)$. The 
\textit{$k$th elementary ideal} $I_k$ of $M$ is the ideal 
$I_k\subset\mathbf{R}$ 
generated by the $(n-k)$-minors of $A$. It is well-known (see for example
\cite{L}) that $I_k$ does not depend on the choice of presentation matrix $A$
for $M$. In particular, the greatest common divisor of the $(n-k)$-minors
of $A$ with respect to a fixed choice of term ordering, denoted
\[\Delta_k(M)=\mathrm{gcd}\{(n-k)\mathrm{-minors\ of\ } A\}\] is a generator 
for the minimal principal ideal $P_k$ containing $I_k$.
If $A$ is a square matrix, then $\Delta_0(M)$ is simply the determinant of $M$.

\begin{definition}\textup{
Let $X$ be a finite birack for rank $N$, $\mathbf{R}$ a Laurent polynomial 
ring and $M$ a $\mathbb{Z}[X]$-module structure on $\mathbf{R}$. The $k$th
\textit{polynomial birack module enhanced invariant} of an oriented link
$L$ is the multiset
\[\Phi_{X}^{M,\Delta_k}(L)=\{\Delta_k(M[f])\ |\ f\in \mathcal{L}(L,X)\}\]
where $\mathcal{L}(L,X)$ is the set of $X$-labelings of diagrams of $L$
over a complete tile of framings mod $N$.
}\end{definition}

We note that $\Phi^{M,\Delta_k}_X$ extends to virtual knots and links in the 
usual way, i.e., by ignoring the virtual crossings. See \cite{K} for more. 
By construction, we have:

\begin{theorem}
If $L$ and $L'$ are ambient isotopic classical or virtual links, $X$ is a 
finite birack, and $M$ is a polynomial birack module, then for each 
$k=0,1,2,\dots$ we have
\[\Phi_{X}^{M,\Delta_k}(L)=\Phi_{X}^{M,\Delta_k}(L').\]
\end{theorem}

\begin{remark}\textup{Note that, like the usual Alexander polynomial,
the polynomials in $\Phi_{X}^{M,\Delta_k}$ are only defined up to multiplication
by units in $\mathbf{R}$. We can normalize these polynomials by multiplying by
units to obtain a polynomial with constant term $1$ for ease
of comparison. For instance, if $\mathbf{R}=\mathbb{Z}_5[q^{\pm 1}]$, then the 
polynomial $4q^2+3q^3+q^4$ would normalize to $1+2q+4q^2$ after multiplication 
by $4q^{-2}$.
}\end{remark}

\begin{example}\textup{Let $X=\{1\}$ be the birack of one element and
$\mathbf{R}=\mathbb{Z}[t^{\pm 1},r^{\pm 1}]$. Then $X$ has a birack module 
structure given by the matrix 
\[M_{\mathbf{R}}=\left[\begin{array}{r|r|r} t & 1-tr & r \end{array}\right].\] 
The 0th polynomial birack module enhanced invariant $\Phi_X^{M,\Delta_0}$ is then
a singleton set whose entry is the \textit{generalized Alexander polynomial},
also known as the \textit{Sawollek polynomial} after a different normalization;
see \cite{KR,S} for more. 
}
\end{example}

\begin{example}\textup{Let $X=\{1\}$ be the birack of one element and
$\mathbf{R}=\mathbb{Z}[t^{\pm 1}]$. The $\mathbf{R}$-module structure with 
matrix 
\[M_{\mathbf{R}}=\left[\begin{array}{r|r|r} t & 1-t & 1 \end{array}\right]\]
has $k=1$ polynomial birack module enhanced invariant $\Phi_X^{M,\Delta_1}$ 
a singleton set whose entry is the usual Alexander polynomial $\Delta(K)$. 
Indeed, $\Phi_X^{M,\Delta_k}$ in this case has single entry given by the usual 
$k$th Alexander polynomial $\Delta_k(K)$.
}\end{example}

Thus, we can think of a polynomial birack module structure on a ring $R$
as determining a customized Alexander polynomial for each $X$-labeling of
our link $L$, and the invariant $\Phi_{X}^{M,\Delta_k}$ collects these polynomials
to form a new invariant whose cardinality is the integral birack counting
invariant but whose entries can distinguish knots and link which have the
same integral birack counting invariant.

\section{\large\textbf{Computations and Applications}}\label{CA}

In this section we collect a few examples of the new invariant. We begin with 
an illustration of how the invariant is computed.

\begin{example}\label{ex1}\textup{
Let $X$ be the birack with birack matrix
\[X=\left[\begin{array}{rr|rr}
1 & 1 & 2 & 2 \\
2 & 2 & 1 & 1 \\
\end{array}\right]\]
and let $\mathbf{R}=\mathbb{Z}_5[q^{\pm 1}]$. A computer search using our 
custom \texttt{python} code found birack modules over $\mathbf{R}$ including 
the module specified by the matrix
\[M_R=\left[\begin{array}{rr|rr|rr}
q & q & 1+2q & 2+4q & 1 & 1 \\
q & q & 3+q & 1+2q & 4 & 4
\end{array}\right].\] 
The birack rank of $X$ is $2$, so for any link $L$ of $c$ components,
we must find $X$-labelings for all writhe vectors in $(\mathbb{Z}_2)^c$.
The virtual trefoil, denoted $2.1$ in the virtual knot table on the 
knot atlas \cite{KA}, has two $X$-labelings $f_1,f_2$ for even writhe 
and no $X$-labelings for odd writhe. These labelings determine the listed
presentation matrices for the modules of bead labelings:
\[\raisebox{-0.5in}{\includegraphics{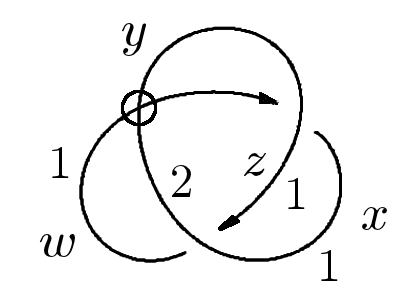}}
\quad \Rightarrow\quad
\left[\begin{array}{rrrr}
s_{11} & 0 & t_{11} & -1 \\
r_{11} & -1 & 0 & 0 \\
-1 & s_{21} & 0 & t_{21} \\
0 & r_{21} & 0 & -1 \\ 
\end{array}\right]
\quad\Rightarrow\quad M_{f_1}=
\left[\begin{array}{cccc}
1+2q & 0 & q & 4 \\
1 & 4 & 0 & 0 \\
4 & 3+q & 0 & q \\
0 & 4 & 0 & 4 \\ 
\end{array}\right]
\]
\[\raisebox{-0.5in}{\includegraphics{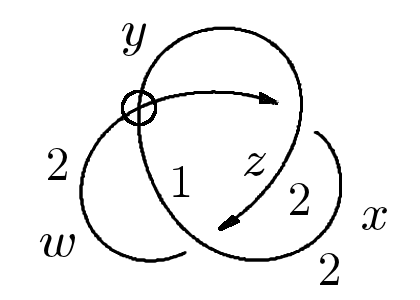}} 
\quad \Rightarrow\quad
\left[\begin{array}{rrrr}
s_{22} & 0 & t_{22} & -1 \\
r_{22} & -1 & 0 & 0 \\
-1 & s_{12} & 0 & t_{12} \\
0 & r_{12} & 0 & 4 \\ 
\end{array}\right]
\quad\Rightarrow\quad M_{f_2}=
\left[\begin{array}{cccc}
1+2q & 0 & q & 4 \\
4 & 4 & 0 & 0 \\
4 & 2+4q & 0 & q \\
0 & 1 & 0 & 4 \\ 
\end{array}\right]
\]
so the $k=0$ invariant after normalization is 
\[\Phi_{X}^{M,\Delta_0}=\{\det(M_{f_1}),\det(M_{f_2})\}=
\{1+q+3q^2,1+q+3q^2\}=\{2\times(1+q+3q^2)\}.\]
Note that since $\Phi_{X}^{M,\Delta_0}(0_1)=\{2\times 0\}$, this example shows that
$\Phi_X^{M,\Delta_0}$ is not determined by the integral birack counting invariant
and hence is a proper enhancement.
}\end{example}

Our next example demonstrates that $\Phi_X^{M,\Delta_0}$ is not determined by the
generalized Alexander polynomial by distinguishing two virtual knots which
have the same generalized Alexander polynomial.

\begin{example}\textup{Let $X, R$ and $M$ be as in example \ref{ex1}. Then
the virtual knots numbered $4.10$ and $4.17$ in the knot atlas \cite{KA}
both have generalized Alexander polynomial $(1 + r)(1 - r)(1 - t)(1 - rt)$, 
but are distinguished by $\Phi_{X}^{M,\Delta_0}$:
\[\begin{array}{cc}
\includegraphics{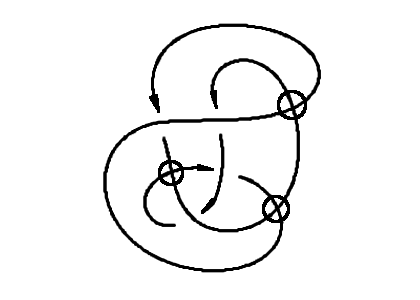} &
\includegraphics{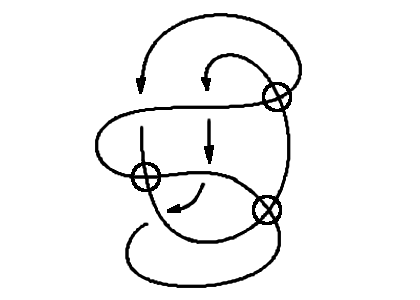} \\
\Phi_{X}^{M,\Delta_0}(4.10) = \{ 2\times (1+2q+4q^2+3q^3) \}
&
\Phi_{X}^{M,\Delta_0}(4.17) = \{2\times (1+q+3q^2)\} \\
\end{array}\]
}\end{example}

Our final example gives an impression of the effectiveness of 
$\Phi_x^{M,\Delta_k}$ as an invariant by sampling the ability of 
$\Phi_x^{M,\Delta_k}$ to differentiate knots for a randomly selected birack module.

\begin{example}\textup{Let $X$ be the birack in example \ref{ex1}.
We randomly selected a polynomial birack $X$-module over 
$R=\mathbb{Z}[q^{\pm 1}]$, given by the matrix
\[M_R=\left[\begin{array}{rr|rr|rr}
q & q & 1+q & 3+3q & 2 & 2 \\
q & q & 2+2q & 1+q & 3 & 3
\end{array}\right]\]
and computed $\Phi_{X}^{M,\Delta_0}$ for the virtual knots with 4 and fewer
crossings as listed on the knot atlas (\cite{KA}) using our 
\texttt{python} code, available at \texttt{www.esotericka.org}. The 
results are collected in the following table.
\[
\begin{array}{r|l}
\Phi_{X}^{M,\Delta_0}(L) & L \\ \hline
\{2\times (0)\} & 3.1,3.5,3.6,3.7,4.8,4.10,4.16,4.32,4.41,4.47,4.50,4.55,4.56,4.58,4.59,\\
& 4.68,4.70,4.71,4.72,4.75,4.76,4.77,4.85,4.86,4.89,4.90,4,96,4.98,4.99,\\
 & 4.102,4.105,4.106,4.107,4.108 \\
\{2\times (1+q+2q^2+4q^3+2q^4)\} & 4.29 \\
\{2\times (1+q+4q^2+4q^3\} & 4.6,4.13,4.17,4.19,4.23,4.24,4.26,4.31,4.35,4.42,4.46,4.51,4.57,4.66, \\ 
& 4.67,4.79,4.93,4.97,4.103 \\
\{2\times (1+q+4q^3+4q^4)\} & 4.9,4.15,4.37,4.45,4.69,4.78,4.92,4.95,4.104 \\
\{2\times (1+2q+4q^2+3q^3)\} & 4.74 \\
\{2\times (1+2q+3q^3+4q^4)\} & 4.12,4.21,4.36,4.61,4.65,4.73 \\
\{2\times (1+3q+4q^2+2q^3)\} & 4.30 \\
\{2\times (1+4q+4q^2+q^3)\} & 4.3,4.25 \\
\{2\times (1+4q^2)\} & 3.2,3.3,4.1,4.4,4.5,4.11,4.14,4.18,4.20,4.22,4.27,4.28,4.33,4.34,4.38,\\
& 4.39,4.40,4.43,4.44,4.48,4.49,4.52,4.54,4.60,4.62,4.63,4.64,4.81,4.82, \\
& 4.84,4.87,4.88,4.94,4.101\\
\{2\times (1+3q^2+q^4)\} & 4.2 \\
\{2\times (1+4q^4)\} & 4.7,4.53,4.80,4.91,4.100 \\
\end{array}
\]
}\end{example}

\section{\large\textbf{Questions}}\label{Q}

In this section we collect a few questions and directions for future work.

The Alexander and generalized Alexander polynomial satisfy the well-known
Conway skein relation. For which biracks $X$, modules $M$ and nonnegative 
integers $k$ do the elements of $\Phi_X^{M,\Delta_k}$ satisfy a skein relation? 

The examples we have selected for inclusion in this paper use biracks and 
polynomial birack modules of small cardinality for speed of computation and
convenience of presentation; even so, the $\Phi^{M,\Delta_k}_X$ invariant 
appears to be quite effective at distinguishing virtual knots with only 
these small $M$ and $X$.
We anticipate that biracks of larger cardinality and polynomial birack modules
with more variables or over larger rings should be even more effective
at distinguishing knots and links. Thus, fast algorithms for computing
$\Phi_X^{M,\Delta_k}$ for larger $k$, $M$ and $X$ are of interest.

\bigskip

\noindent
\textsc{Department of Mathematical Sciences \\
Claremont McKenna College \\
850 Columbia Ave. \\
Claremont, CA 91711}

\end{document}